\documentclass[graybox]{svmult}

\usepackage{type1cm}        
\usepackage{makeidx}         
\usepackage{graphicx}        
\usepackage{multicol}        
\usepackage[bottom]{footmisc}

\usepackage{newtxtext}       %
\usepackage{newtxmath}       
\usepackage{wrapfig}

\makeindex             

\usepackage{xcolor}

\usepackage{soul} 
\usepackage{dsfont} 
\usepackage[hidelinks]{hyperref}

\newcommand{\curl}[1]{\text{curl}\left(#1\right)}
\newcommand{\dd}{\mathsf{d}}

\newcommand{\RR}{\mathbb{R}}
\newcommand{\CC}{\mathbb{C}}
\newcommand{\Hcurl}[1]{\mathbf{H}\left(\text{curl, }#1\right)}

\newcommand{\trace}[1]{\gamma^T\left(#1\right)}
\newcommand{\curltrace}[1]{\gamma^t\left(#1\right)}

\newcommand{\tracei}[2]{\gamma^T_#1\left(#2\right)}
\newcommand{\curltracei}[2]{\gamma^t_#1\left(#2\right)}

\begin{document}

\title*{Parallel domain decomposition solvers for the time harmonic Maxwell equations}
\author{Sven Beuchler, Sebastian Kinnewig, Thomas Wick}
\institute{Sven Beuchler, Sebastian Kinnewig, Thomas Wick \at Leibniz University Hannover, Institute of Applied Mathematics, Welfengarten 1, 30167 Hannover, Germany 
\{beuchler,kinnewig,wick\}@ifam.uni-hannover.de
\and and Cluster of Excellence PhoenixD (Photonics, Optics, and
	Engineering - Innovation Across Disciplines), Leibniz Universit\"at Hannover, Germany
}
%
%
\maketitle

\abstract{
  The time harmonic Maxwell equations are of current interest in 
  computational science and applied mathematics with many applications
  in modern physics.
  In this work, we present
  parallel finite element solver for the
  time harmonic Maxwell equations and compare 
  different preconditioners.
  We show numerically that standard preconditioners like incomplete LU and the
  additive Schwarz method lead to slow convergence for iterative solvers like
  generalized minimal residuals, especially for high wave numbers. 
  Even more we show that also more
  specialized methods like the Schur complement method 
  also yield slow convergence.
  As an example for a highly adapted solver for the time harmonic Maxwell 
  equations we use a combination of a block preconditioner and a domain 
  decomposition method (DDM),
  which also preforms well for high wave numbers. 
  Additionally we discuss briefly further approaches to solve high frequency 
  problems more efficiently. 
  Our developments are done in the open-source finite element library deal.II.
}

\section{Introduction}
\label{sec:1}
The time harmonic Maxwell (THM) equations are of great interest in 
applied mathematics \cite{HaKuLaReiSchoe01,bk:monk:03,art:gray:15,art:bou:15,Bue13,LaPauRe19} 
and current physics applications, e.g., the excellence cluster 
PhoenixD.\footnote{\url{https://www.phoenixd.uni-hannover.de/en/}}
However, the numerical solution is challenging. This is specifically true for 
small wave numbers. Various solvers and preconditioners have been proposed, while the most promising 
are based on domain decomposition methods (DDM) \cite{ToWid05}. 
In \cite{art:bou:15}, a quasi-optimal domain decomposition (DD) algorithm was proposed, 
mathematically analyzed and demonstrated to perform well for several numerical examples. 

The goal of this work is to employ  
the domain decomposition method from \cite{art:bou:15}
and to re-implement the algorithm in the modern finite 
element library deal.II \cite{art:dealii:20}.
Therein, the construction of the subdomain interface conditions 
is a crucial aspect for 
which we use Impedance Boundary Conditions.
{
Instead of handling the resulting linear system with a direct solver, which is typically done
for the THM, we apply a well chosen {block preconditioner} to the linear system so we can solve it 
with an iterative solver like GMRES (generalized minimal residuals). 
Additionally high polynomial N\'ed\'elec elements are
used in the implementation of the DDM, see \cite{phd:zag:06}.
}

This implementation is computationally 
compared to several other (classical) preconditioners such 
as incomplete LU, additive Schwarz, Schur complement. 
These comparisons are done for different wave numbers. Higher wave numbers are well-known
to cause challenges for the numerical solution.


The outline of this work is as follows:
In the section \ref{sec:2} we introduce some notation. In section \ref{sec:3} we introduce 
the  domain decomposition method {(DDM)} for the THM, 
furthermore we introduce a 
block preconditioner which will allow us 
to solve the THM with iterative solves instead of direct solvers inside of DDM.
In the section \ref{sec:4} we will compare the results of the block preconditioner 
with the performance of {different} preconditioners. Moreover we will present 
some results of the combination of the preconditioner and the DDM for two benchmark problems.

%

\section{Equations and finite element discretization}
\label{sec:2}
Let $\Omega \subset \RR^d,~ d \in \{2,3\}$ be a bounded domain with 
sufficiently smooth boundary $\Gamma$. 
The latter is partitioned into $\Gamma = \Gamma^{\infty} \cup \Gamma^{\text{inc}}$. 
Furthermore, the time harmonic Maxwell equations are then defined as follows: Find the electric field 
$\vec E \in \Hcurl{\Omega}:=\{v \in \mathcal{L}^2(\Omega), ~\curl{v} \in \mathcal{L}^2(\Omega)\}$ such that
\begin{equation}
    \label{eq:THM:SingleDomain}
    \left\{\begin{array}{l l l}
    \curl{ \mu^{-1} \curl{ \vec E } } - \omega^2 \vec E           & = \vec 0                    & \text{ in } \Omega    \\
    \mu^{-1} \curltrace{\curl{ \vec E }} - i \kappa \omega \trace{\vec E}  & = \vec 0                    & \text{ on } \Gamma^{\infty} \\
    \trace{\vec E}                                                         & = \trace{\vec E^\text{inc}} & \text{ on } \Gamma^{\text{inc}} 
  \end{array}\right. ,
\end{equation} 
where $\vec E^{\text{inc}}:\RR^d \to \CC^d,~ d \in \{2,3\}$
 is some given incident electric 
field, $\omega > 0$ is the wave number which is defined by $\omega:=\frac{2\pi}{\lambda}$, 
where $\lambda > 0$ is the wave length. 
Let $\Omega$ be a domain with smooth interface.

Following \cite{GirRav,bk:monk:03}, we define traces
$\gamma^t:\Hcurl{\Omega} \rightarrow \mathbf{H}_\times^{-1/2}(\operatorname{div}, \Gamma)$ and 
$\gamma^T:\Hcurl{\Omega} \rightarrow \mathbf{H}_\times^{-1/2}(\operatorname{curl}, \Gamma)$ 
by
\begin{equation*}
  \curltrace{\vec v} = \vec n \times \vec v \text{ and } \trace{\vec v} = \vec n \times (\vec v \times \vec n)
\end{equation*}
where the vector $\vec n$ is the normal to $\Omega$,
$\mathbf{H}_\times^{-1/2}(\operatorname{div}, \Gamma) := \{ \vec v \in \mathbf{H}^{-1/2}(\Gamma) : \vec v \cdot \vec n=0, 
\operatorname{div}_\Gamma \vec v \in  \mathbf{H}^{-1/2}(\Gamma) \}$ 
is the space of well-defined surface divergence fields, 
$\mathbf{H}_\times^{-1/2}(\operatorname{curl}, \Gamma) := \{\vec v \in \mathbf{H}_\times^{-1/2} : \vec v \cdot \vec n=0, \operatorname{curl}_\Gamma \vec v\in  \mathbf{H}^{-1/2}(\Gamma) \}$
is the space of well-defined surface curls.

System \eqref{eq:THM:SingleDomain} is called time harmonic, because the time dependence can be 
expressed by $e^{i \omega \tau}$, where $\tau \geq 0$ denotes the time.

For the implementation with the help of a Galerkin finite element method, we need the discrete weak form. 
{Let 
$\mathcal{N}_h^p := \{ v_h \in X : v_h|_K (x) = a_K (x) + \big(x \times b_K(x)\big) , ~ a_K, b_K \in [P^p(K)]^3  
~ \forall~ K \in \tau_h(\Omega)\}$ 
be the N\'ed\'elec space \cite{bk:monk:03}.
Based on the de-Rham cohomology, basis functions can be developed, \cite{phd:zag:06}.}
Find $\vec E_h \in \mathcal{N}_h^p(\Omega)$ such that
\begin{align}
  \label{eq:THM:WeakSingleDomain}
  &\quad \int_\Omega \left( \mu^{-1}\curl{\vec E_h}\curl{\vec \phi_h}  - \omega^2 \vec E_h \vec \phi_h \right) \dd x \nonumber \\
  &+\int_{\Gamma^\infty} i ~\mu^{-1} \trace{\vec E_h} \trace{\vec\phi_h} \dd s 
  = \int_{\Gamma^\text{inc}} \trace{\vec E_h^\text{inc}} \trace{\vec\phi_h} \dd s
\quad\forall \vec \phi \in \mathcal{N}_h^p(\Omega).
\end{align}
In order to obtain a block system for the numerical solution process,
we define the following elementary integrals
\begin{equation}
\begin{array}{l c l}
  (A)_{u,v} & = & \int_\Omega \mu^{-1}\curl{\vec \phi_u}\curl{\vec \phi_v}  \\
  (M)_{u,v} & = & \int_\Omega \vec \phi_u \vec \phi_v  \\
  (B)_{u,v} & = & \int_{\Gamma^\infty} i \mu^{-1}\trace{\vec \phi_u}\trace{\vec \phi_v} \\
  (s)_{u}   & = & \int_{\Gamma^\text{inc}} \trace{\vec E^\text{inc}}\trace{\vec \phi_u}, \\
\end{array}
\end{equation}
where $\phi_u,~\phi_v \in \mathcal{N}_h^p(\Omega)$. 
To this end, System \eqref{eq:THM:SingleDomain} can be written in the form
\begin{equation}
  \label{eq:THM:Block}
  \begin{pmatrix}
    A - \omega^2 M      & -B \\
    B                   & A - \omega^2 M
  \end{pmatrix}
  \begin{pmatrix}
    \vec E_{RE} \\
    \vec E_{IM}
  \end{pmatrix}
  =
  \begin{pmatrix}
    \vec s_{RE} \\
    \vec s_{IM}
  \end{pmatrix},
\end{equation}
where $\vec E = \vec E_{RE} + i \vec E_{IM} $ and $\vec s = \vec s_{RE} + i \vec s_{IM}$, where $i$ denotes the imaginary number.

\section{Numerical solution with domain decomposition and preconditioners}
\label{sec:3}

\subsection{Domain decomposition}
Due to the ill-posed nature of the time harmonic Maxwell equations,
{a  successful approach to solve the THM is based on the DDM \cite{ToWid05}.
As the name suggests, the domain is divided into smaller subdomains. 
As these subdomains become small enough they can be handled by a direct solver.}
To this end, we divide the domain as follows:
\begin{equation}
  \Omega = \bigcup_{i=0}^{N_{\text{dom}}} \Omega_i
\end{equation}
where $N_{\text{dom}}$ is the number of domains, 
since we consider a non-overlapping DDM 
$\Omega_i \cap \Omega_j = \emptyset$, if $i \neq j ~ \forall i,j \in \{1,\ldots, N_\text{dom}\}$ and
we denote the interface from two neighbouring cells by
$\partial \Omega_i \cap \partial \Omega_j = \Sigma_{ij} = \Sigma_{ji}, ~ \forall i,j \in \{1,\ldots, N_\text{dom}\}$.

The second step of the DD is an iterative method, indexed by $k$,
to compute the overall electric field $\vec E$.
Therefore we begin by solving System \eqref{eq:THM:SingleDomain} 
on each subdomain 
$\Omega_i$, we denote the solution of every subsystem by $\vec E^{k = 0}_{i}$.
From this we can compute the first interface condition by
\begin{equation}
  \label{eq:interface}
  g^{k = 0}_{ji} := -\mu^{-1} \curltracei{i}{\curl{E_i^{k = 0}}} - i k S\left(\tracei{i}{E_i^{k = 0}}\right),
\end{equation}
where $S$ describes some boundary operator, 
which we will discuss in more detail below.
Afterward, we obtain the next iteration step $\vec E_i^{k+1}$ via:
\begin{eqnarray}
  \label{eq:DD:iteration}
  \left\{\begin{array}{l l l}
    \curl{ \mu^{-1} \curl{ \vec E^{k+1}_i } } - \omega^2 \vec   E^{k+1}_i & = \vec 0                      & \text{ in } \Omega_i   \\[2mm]
    \mu^{-1} \curltracei{i}{ \curl{ \vec E^{k+1}_i } } - i \kappa \omega \tracei{i}{ \vec E^{k+1}_i }  & = 0        & \text{ on } \Gamma^{\infty}_i \\[2mm]
    \tracei{i}{ \vec E^{k+1}_i }                                                                  & = \vec E^\text{inc}_i & \text{ on } \Gamma^{\text{inc}}_i \\[2mm]
    \mu^{-1} S\left( \curltracei{i}{\curl{ \vec E^{k+1}_i }} \right) - i \kappa \omega \tracei{i}{ \vec E^{k+1}_i }   & = g^k_{ji}        & \text{ on } \Sigma_{i,j} 
  \end{array}\right.
\end{eqnarray}
Once $\vec E^{k+1}_i$ is computed, the interface is updated by
\begin{equation}
  \label{eq:interface_update}
  g_{ji}^{k+1} = -\mu^{-1} \curltracei{i}{\curl{E_i^{k+1}}} - i k S\left(\tracei{i}{E_i^{k+1}}\right) 
               = -g_{ij}^k - 2 i k S\left(\tracei{i}{E_i^{k+1}}\right) 
\end{equation}
where $\vec E^k_i \rightarrow \vec E|_{\Omega_i}$ as $k \rightarrow \infty$.
This convergence depends strongly on the chosen surface operator $S$. 
For a convergence analysis when the IBC are considered, see 
\cite{art:dol09}

This iteration above can be interpreted as one step of the Jacobi fixed point method for the linear system
\begin{equation}
  \label{eq:iteration}
  (\mathds{1} - \mathcal{A}) \vec g = \vec b
\end{equation}
where $\mathds{1}$ is the identity operator, $\vec b$ is the vector of the incident electric field,
$\mathcal{A}$ is defined by $\mathcal{A}\vec g^k = \vec g^{k+1}$ 
and Equations \eqref{eq:DD:iteration}, \eqref{eq:interface_update}.
Convergence is achieved for
$\|(1 - \mathcal{A})\vec g^k - b\| < TOL$ with some small tolerance $TOL>0$. 
Often, $TOL = 10^{-6},\ldots, 10^{-8}$.
Instead of a Jacobi fixed point method one can also use a GMRES method to solve 
\eqref{eq:iteration} more efficiently.

The crucial point of the DD is the choice of the interface conditions between the subdomains. 
The easiest choice is a non-overlapping Schwarz decomposition, where Dirichlet like interface 
conditions are used. 
{For large wave numbers, e.g. the parameter $\omega$ becomes large,
the ellipicity of the system goes lost. Consequently, a convergence of this algorithm for the time harmonic 
Maxwell equations for all $\omega$ cannot be expected; see \cite{art:ern:11,GraSpeVai}.
An analysis for an overlapping additive Schwarz method is given in \cite{Grahametal2019}.}



Rather, we need more sophisticated tools in
which the easiest choice are Impedance Boundary Conditions (IBC), which can be
classified as Robin like interface conditions 
\begin{equation}
  S = \mathds{1}.
\end{equation}

\subsection{Preconditioner}

As it is clear, the DDM is an iterative method, where we have to 
solve system \eqref{eq:DD:iteration} on each 
subdomain in each iteration step $k$. 
Usually, this is done by 
a direct solver, but instead, we can 
use a GMRES solver, which is preconditioned by an 
approximation of the block system
\begin{equation}
  \label{eq:preconditioner}
  \begin{pmatrix}
    A - \omega^2 M & 0 \\
    0 & A - \omega^2 M
  \end{pmatrix}^{-1}.
\end{equation}
Therefore we need to compute an approximation 
of $(A - \omega^2 M)^{-1}$, and we obtain this approximation
by applying the AMG preconditioner provided by 
MueLu \cite{MueLu}, where for the level transitions a
direct solver is used. The latter is necessary, since 
otherwise the AMG preconditioner does not perform well for the THM.
On the one hand, this procedure is cost expensive. On the other hand,
we can reuse the preconditioner each time we solve 
system \eqref{eq:DD:iteration}.

An other possible choice is to use the AMG preconditioner to compute
directly an approximation of
\begin{equation*}
  \begin{pmatrix}
    A - \omega^2 M & -B \\
    -B & A - \omega^2 M
  \end{pmatrix}^{-1}.
\end{equation*}
With this preconditioner only a few GMRES iterations are needed to solve the system 
\eqref{eq:DD:iteration}. Since we computing an approximation of the complete inverse
this comes with much higher memory consumption, than using \eqref{eq:preconditioner} 
as preconditioner. Actually the memory consumption while using an iterative solver 
with \eqref{eq:preconditioner} as an preconditioner is even lower, than the memory
consumption from a direct solver, which we show numerically in the next chapter.
Therefore the block diagonal preconditioner is used in the following.


\section{Numerical tests}
\label{sec:4}
In this section, we compare the performance of different preconditioners 
for two numerical examples. We choose a simple wave guide 
as our benchmark problem, 
moreover we test the performance of our method on a Y beam splitter. 
Our implementation is based on 
the open-source finite element library 
\texttt{deal.II} \cite{art:dealii:20} with Trilinos \cite{trilinos1} and MueLu \cite{MueLu}. As a direct solver, MUMPS (Multifrontal Massively Parallel Sparse Direct Solver) \cite{MUMPS} is used.
We perform an additive domain decomposition and compute each step in parallel with MPI.
For the computations an Intel Xeon Platinum 8268 CPU was used with up to 32 cores.

\subsection{Example 1: Block benchmark}
\label{sec:4.1}

\begin{wrapfigure}{r}{4cm}
  \includegraphics[scale=.34]{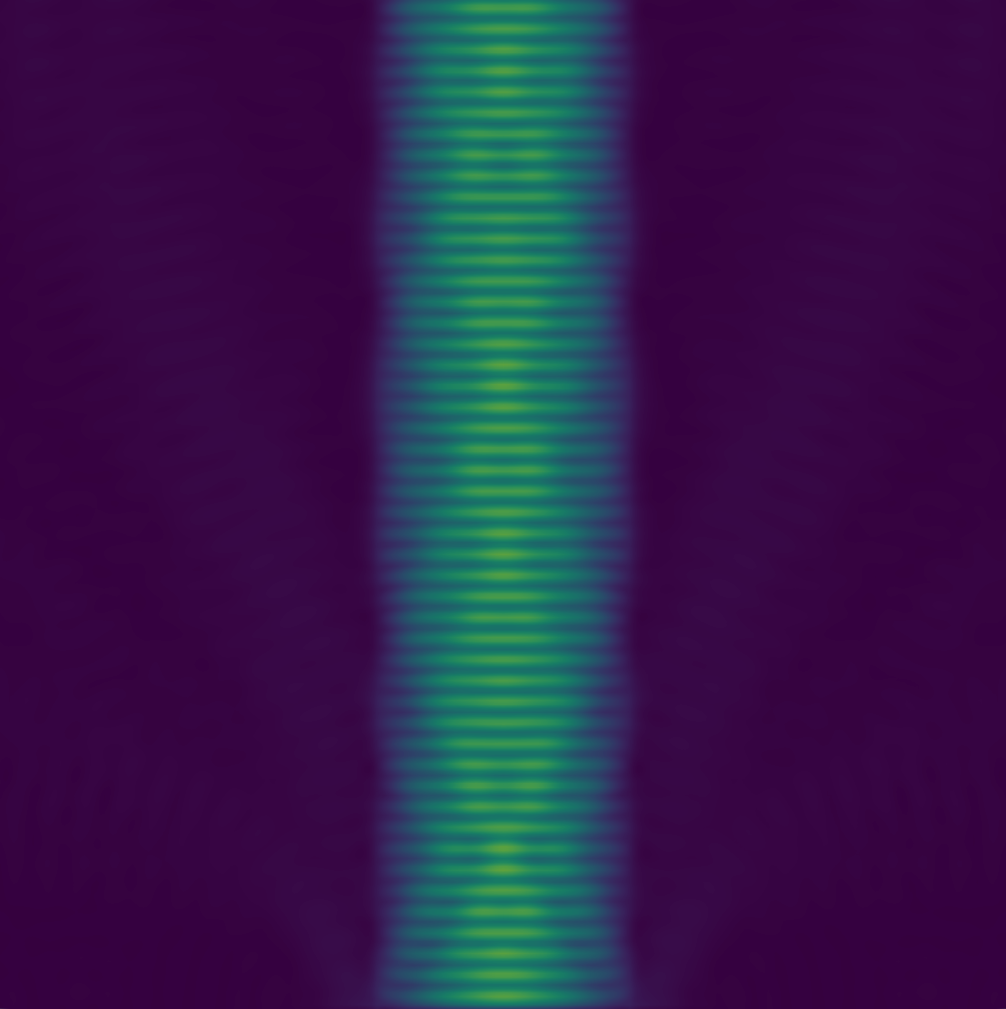}
  \caption[t]{{
    As a benchmark problem, we consider a $1 \times 1$ square with different wave numbers 
    (here $\lambda = 50$). In the center is a material with the refractive index $n_\text{center} = 1.516$
    and as cladding the refractive index of air was used $n_\text{air} = 1.0$.
  }}
  \label{fig:NumTest:block}       
\end{wrapfigure}

Let us consider a squared domain decomposed of a material with a higher refractive index in the center a
carrier material with a lower refractive index beside it, see Figure \ref{fig:NumTest:block}.

Table \ref{tab:NumTest:iter} displays the GMRES iterations with a relative accuracy of $\epsilon = 10^{-8}$
for differenent preconditioners:
\begin{itemize}
\item ILU, incomplete LU decomposed of \eqref{eq:THM:Block},
\item the implemented additive Schwarz preconditioner of \cite{art:dealii:20, trilinos1} 
  \footnote{\url{https://www.dealii.org/current/doxygen/deal.II/classTrilinosWrappers_1_1PreconditionSSOR.html}},
\item a Schur complement preconditioner based on 
  \footnote{\url{https://www.dealii.org/current/doxygen/deal.II/step_22.html}},
\item the block preconditioner \eqref{eq:preconditioner}.
\end{itemize}
Except for the last case, the GMRES iteration numbers grow for large $\omega$.

\begin{table}[h!]
  \caption{
    {Example 1: GMRES iterations with different preconditioners.} 
  }
  \label{tab:NumTest:iter}   

  \begin{tabular}{p{2.2cm}p{1.8cm}p{2.2cm}p{2.4cm}p{2.6cm}} 
    \hline\noalign{\smallskip}
    wave number $\omega$ & \multicolumn{4}{l}{GMRES iterations with the preconditioner} \\ 
                         & ILU    & additive Schwarz & Schur complement & block preconditioner \\ 
    \noalign{\smallskip}\svhline\noalign{\smallskip}
       5.0               & 369     & 34              & 89    &  43 \\ 
      10.0               & 858     & 113             & 166   &  49 \\ 
      20.0               & >1000   & >1000           & 493   &  33 \\ 
      40.0               & -       & -               & 598   &  32 \\ 
      80.0               & -       & -               & >1000 &  28 \\ 
  \noalign{\smallskip}\hline\noalign{\smallskip}
  \end{tabular}
\end{table}

\begin{figure}[h!]
  \label{fig:NumTest:runtime}       
  \begin{minipage}{0.49\textwidth}
    \includegraphics[width=1.0\textwidth]{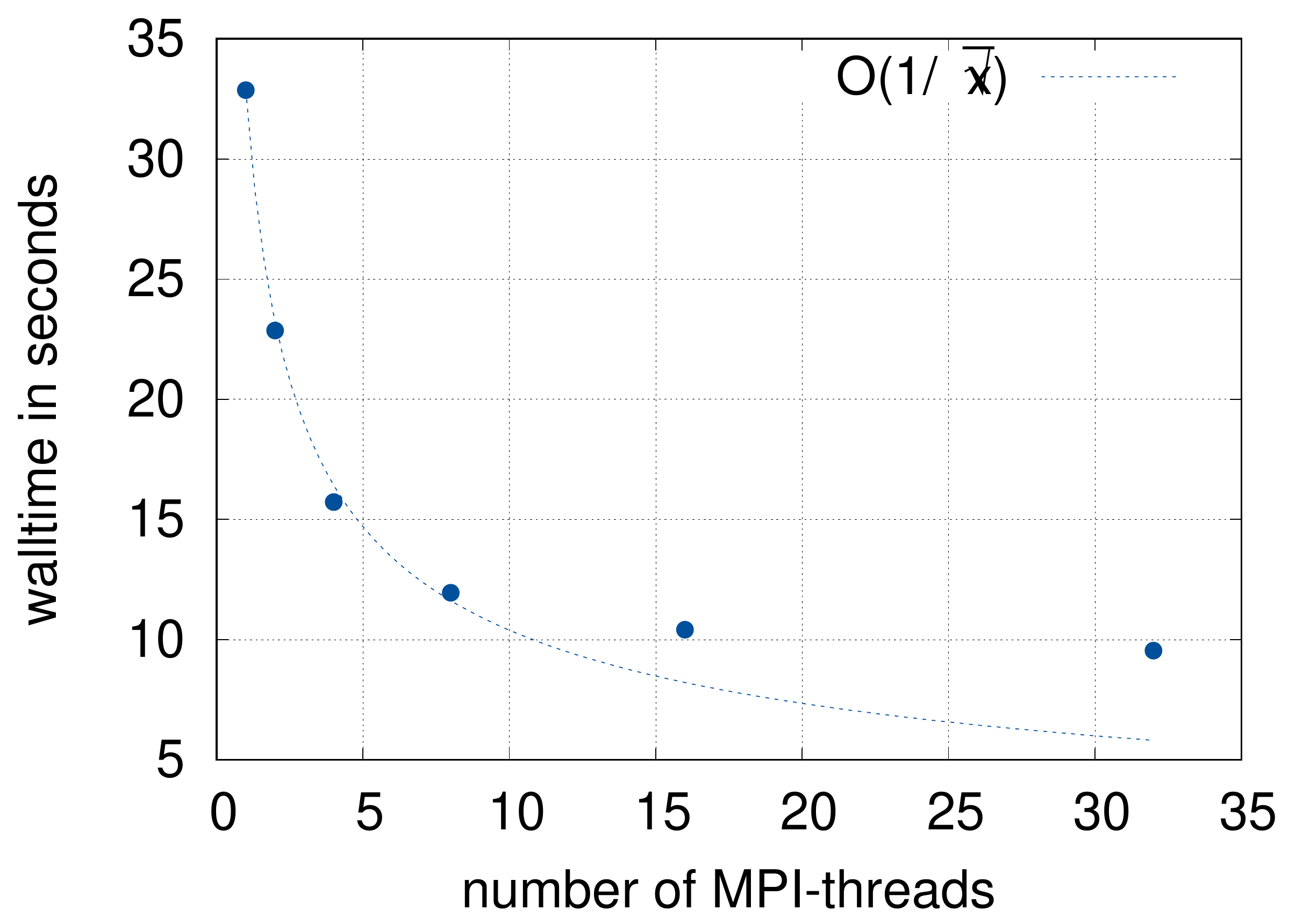}\\
  \end{minipage}
  \hfill
  \begin{minipage}{0.49\textwidth}
    \includegraphics[width=1.0\textwidth]{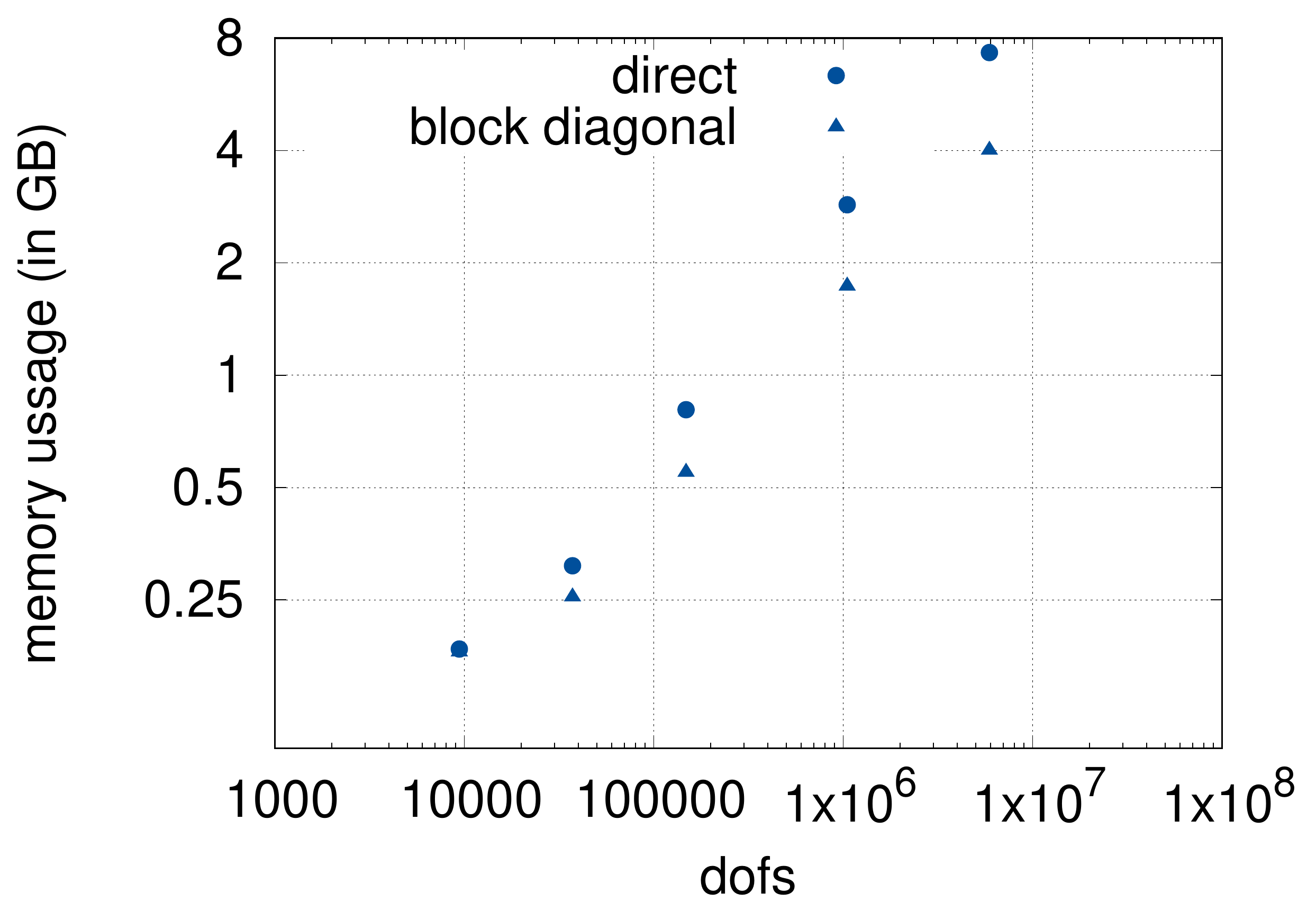}\\
  \end{minipage}
  \ \\
  \ \\
  \caption{{On the left side: memory usage in dependence of the number of dofs. 
  On the right side: walltime in dependence of the number of MPI-threads}.}
\end{figure}

\newpage
\subsection{Example 2: Y beam splitter}
\label{sec:4.2}
Similar as in the simple wave guide, we consider for the Y beam splitter an 
material with a 
higher refractive index placed inside of an carrier material with a 
lower refractive index.
Here we consider a 3D model of a Y beam splitter. 
The mesh was divided into 9 subdomains, 
and the average number of GMRES iterations to solve the 
subdomains are given in table \ref{tab:NumTest:ddm_iter},
are for the wave number $\lambda = 20$.

\begin{table}[h!]
  \caption{{Example 2:} GMRES iterations on each domain}
  \label{tab:NumTest:ddm_iter}   

  \begin{tabular}{l p{0.7cm} p{0.7cm} p{0.7cm} p{0.7cm} p{0.7cm} p{0.7cm} p{0.7cm} p{0.7cm} p{0.7cm}}
    \hline\noalign{\smallskip}
    subdomain id & 1 & 2 & 3 & 4 & 5 & 6 & 7 & 8 & 9 \\
    \noalign{\smallskip}\svhline\noalign{\smallskip}
    average number of GMRES iterations &
    34 & 40 & 41 & 31 & 35 & 39 & 37 & 33 & 32 \\
    \noalign{\smallskip}\hline\noalign{\smallskip}
  \end{tabular}
\end{table}

\begin{figure}[h!]
  \begin{minipage}{0.79\textwidth}
    \includegraphics[width=1.0\textwidth]{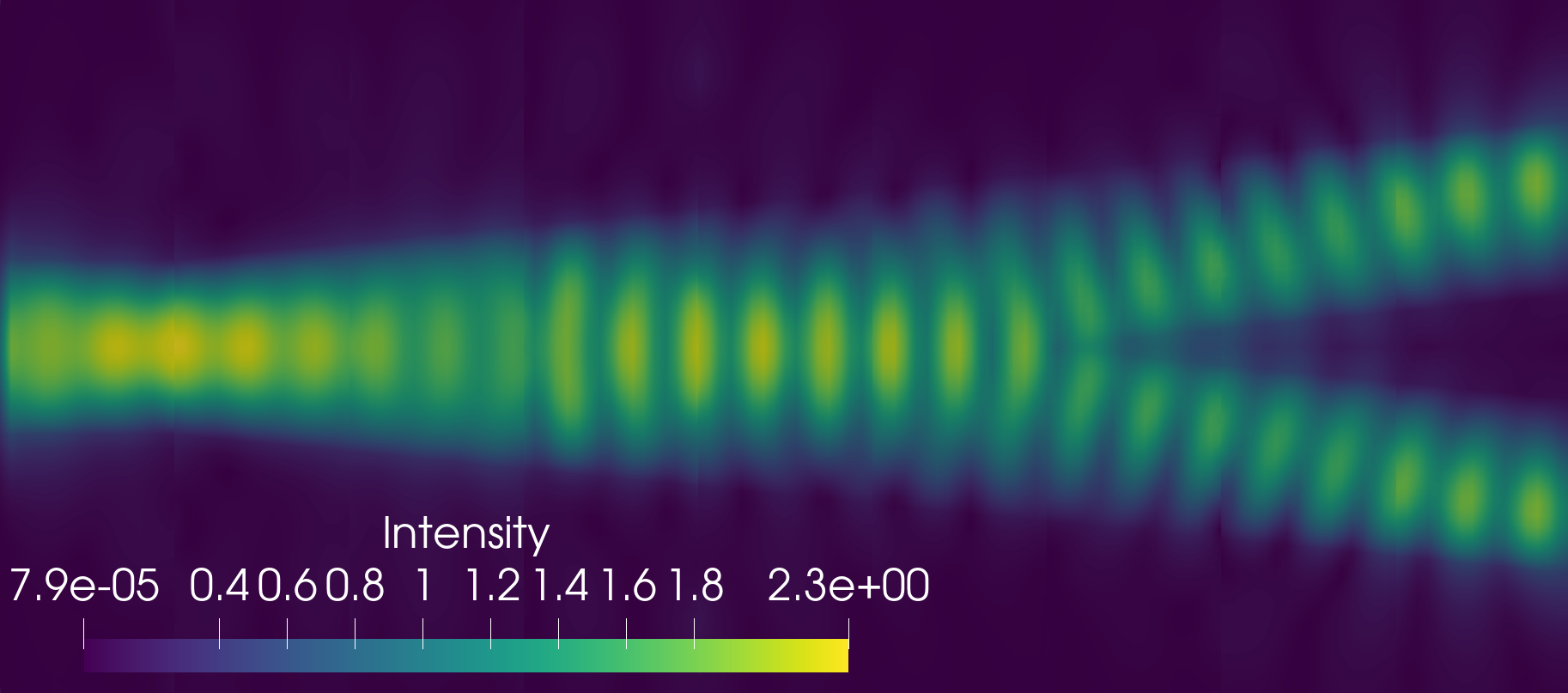}
  \end{minipage}
  \hfill
  \begin{minipage}{0.2\textwidth}
    \includegraphics[width=1.75\textwidth, angle=-90]{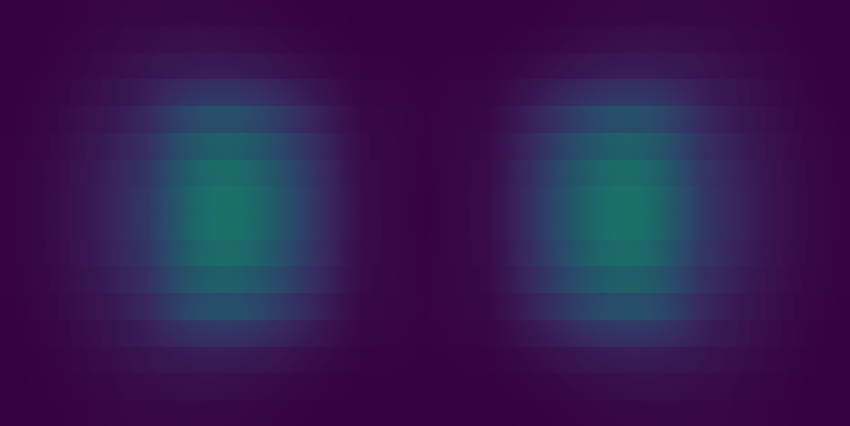}
  \end{minipage}
  \caption{{Intensity plot of the y beam splitter, on the left side is the intensity on the x-y plane and on the right side is the intensity at the output.}}
\end{figure}


\section{Conclusion}
\label{sec:5}
In this contribution, we implemented a domain decomposition 
method with a block preconditioner for the time harmonic Maxwell equations.
Therein, a crucial aspect is the construction of the subdomain 
interface conditions. Our algorithmic developments are demonstrated 
for two configurations of practical relevance, namely a block benchmark 
and a Y beam splitter.

\begin{acknowledgement}
  Funded by the Deutsche Forschungsgemeinschaft (DFG) under Germany’s Excellence Strategy within 
  the Cluster of Excellence PhoenixD (EXC 2122, Project ID 390833453).\\
\end{acknowledgement}


 \bibliographystyle{abbrv}
 \bibliography{lit}

\begin{thebibliography}{10}

\bibitem{MUMPS}
P.~R. Amestoy, A.~Buttari, J.-Y. L'Excellent, and T.~Mary.
\newblock Performance and scalability of the block low-rank multifrontal
  factorization on multicore architectures.
\newblock {\em ACM Trans. Math. Software}, 45(1):Art. 2, 26, 2019.

\bibitem{art:dealii:20}
D.~Arndt, W.~Bangerth, B.~Blais, and et~al.
\newblock The deal.{II} library, {V}ersion 9.2.
\newblock {\em J. Numer. Math.}, 28(3):131--146, 2020.

\bibitem{MueLu}
L.~Berger-Vergiat, C.~A. Glusa, J.~J. Hu, M.~Mayr, A.~Prokopenko, C.~M.
  Siefert, R.~S. Tuminaro, and T.~A. Wiesner.
\newblock {M}ue{L}u multigrid framework.
\newblock \url{http://trilinos.org/packages/muelu}, 2019.

\bibitem{Grahametal2019}
M.~Bonazzoli, V.~Dolean, I.~G. Graham, E.~A. Spence, and P.-H. Tournier.
\newblock Domain decomposition preconditioning for the high-frequency
  time-harmonic {M}axwell equations with absorption.
\newblock {\em Math. Comp.}, 88(320):2559--2604, 2019.

\bibitem{art:bou:15}
M.~E. Bouajaj, B.~Thierry, X.~Antoine, and C.~Geuzaine.
\newblock A quasi-optimal domain decomposition algorithm for the time-harmonic
  maxwell’s equations.
\newblock {\em Journal of Computational Physics}, 2015.

\bibitem{Bue13}
M.~B{\"u}rg.
\newblock Convergence of an automatic hp-adaptive finite element strategy for
  maxwell's equations.
\newblock {\em Applied Numerical Mathematics}, 72:188--204, 10 2013.

\bibitem{art:dol09}
V.~Dolean, M.~J. Gander, and L.~Gerardo-Giorda.
\newblock Optimized {S}chwarz methods for {M}axwell's equations.
\newblock {\em SIAM J. Sci. Comput.}, 31(3):2193--2213, 2009.

\bibitem{art:ern:11}
O.~G. Ernst and M.~J. Gander.
\newblock Why it is difficult to solve {H}elmholtz problems with classical
  iterative methods.
\newblock 83:325--363, 2012.

\bibitem{GirRav}
V.~Girault and P.-A. Raviart.
\newblock {\em Finite element methods for {N}avier-{S}tokes equations},
  volume~5 of {\em Springer Series in Computational Mathematics}.
\newblock Springer-Verlag, Berlin, 1986.
\newblock Theory and algorithms.

\bibitem{GraSpeVai}
I.~G. Graham, E.~A. Spence, and E.~Vainikko.
\newblock Domain decomposition preconditioning for high-frequency {H}elmholtz
  problems with absorption.
\newblock {\em Math. Comp.}, 86(307):2089--2127, 2017.

\bibitem{art:gray:15}
A.~Grayver and T.~Kolev.
\newblock Large-scale 3d geoelectromagnetic modeling using parallel adaptive
  high-order finite element method.
\newblock {\em GEOPHYSICS}, 80:E277--E291, 11 2015.

\bibitem{HaKuLaReiSchoe01}
G.~Haase, M.~Kuhn, U.~Langer, S.~Reitzinger, and J.~Schöberl.
\newblock {\em Parallel Maxwell Solvers}, pages 71--78.
\newblock Lecture Notes in Computational Science and Engineering, Vol. 18.
  Springer, Germany, 2001.

\bibitem{trilinos1}
M.~A. Heroux, R.~A. Bartlett, V.~E. Howle, R.~J. Hoekstra, J.~J. Hu, T.~G.
  Kolda, R.~B. Lehoucq, K.~R. Long, R.~P. Pawlowski, E.~T. Phipps, A.~G.
  Salinger, H.~K. Thornquist, R.~S. Tuminaro, J.~M. Willenbring, A.~Williams,
  and K.~S. Stanley.
\newblock An overview of the trilinos project.
\newblock {\em ACM Trans. Math. Softw.}, 31(3):397--423, 2005.

\bibitem{LaPauRe19}
U.~Langer, D.~Pauly, and S.~Repin, editors.
\newblock {\em Maxwell's Equations}.
\newblock volume 24 of Radon Series on Computational and Applied Mathematics,
  Berlin. de Gruyter, 2019.

\bibitem{bk:monk:03}
P.~{Monk}.
\newblock {\em Finite Element Methods for Maxwell's Equations}.
\newblock Oxford Science Publications, 2003.

\bibitem{ToWid05}
A.~Toselli and O.~Widlund.
\newblock {\em Domain decomposition methods - algorithms and theory}.
\newblock Volume 34 of Springer Series in Computational Mathematics. Springer,
  Berlin, Heidelberg, 2005.

\bibitem{phd:zag:06}
S.~{Zaglmayr}.
\newblock {\em High Order Finite Element Methods for Electromagnetic Field
  Computation}.
\newblock PhD thesis, Johannes Kepler University Linz, 2006.

\end{thebibliography}

\end{document}